\documentclass{article}[11pt]
\usepackage[applemac]{inputenc}	

\usepackage[pdftex]{graphicx}
\usepackage{subfig}
\usepackage{color}
\usepackage{amsfonts}
\usepackage{amsmath,amssymb,amsthm}

\usepackage{mathrsfs}

\usepackage{geometry}         
\definecolor{darkmidnightblue}{rgb}{0.0, 0.2, 0.4}
\definecolor{coolblack}{rgb}{0.0, 0.18, 0.39}
\definecolor{persianindigo}{rgb}{0.2, 0.07, 0.48}
\definecolor{st.patrick}{rgb}{0.14, 0.16, 0.48}
\definecolor{royalblue}{rgb}{0.0, 0.14, 0.4}
\RequirePackage[colorlinks,linkcolor=royalblue,citecolor=darkmidnightblue,urlcolor=coolblack]{hyperref}

\theoremstyle{plain}

\newtheorem{theo}{Theorem}[section] 
 
\newtheorem{lem}[theo]{Lemma}
\newtheorem{prop}[theo]{Proposition}

\newtheorem{rema}[theo]{Remark}

\def\eps{\varepsilon}
\def\mbf{\mathbf}

\newcommand{\iv}{\text{inv}_2}
\newcommand{\un}{\mathbf{1}}
\newcommand{\C}{\mathbb{C}}

\newcommand{\Bc}{\mathcal{B}}

\newcommand{\E}{\mathbf{E}}
\newcommand{\Ec}{\mathcal{E}}
\newcommand{\Lc}{\mathcal{L}}
\newcommand{\Ic}{\mathcal{I}}

\newcommand{\R}{\mathbb{R}}
\newcommand{\N}{\mathbb{N}}

\newcommand{\Sy}{\mathfrak{S}}

\newcommand{\Rc}{\mathcal{R}}

\newcommand{\Mc}{\mathcal{M}}

\newcommand{\Uc}{\mathcal{U}}

\newcommand{\ts}{\otimes} 

\newcommand{\Wg}{\text{Wg}}

\newcommand{\Sym}{\mathrm{Sym}}

\newcommand{\vp}{\varphi}

\newcommand{\Hc}{\mathcal{H}}
\newcommand{\Dc}{\mathcal{D}}

\newcommand{\Pc}{\mathcal{P}}

\newcommand{\esp}{\mathbb{E}}

\newcommand{\SO}{\mathrm{SO}}
\newcommand{\OO}{\mathrm{O}}

\newcommand{\SU}{\mathrm{SU}}
\newcommand{\U}{\mathrm{U}}

\newcommand{\Sp}{\mathrm{Sp}}

\newcommand{\usp}{{\mathfrak{sp}(N)}}

\newcommand{\oN}{{\mathfrak{o}(N)}}
\newcommand{\uN}{{\mathfrak{u}(N)}}

\newcommand{\glN}{{\mathfrak{gl}_N}}

\newcommand{\su}{\mathfrak{su}}

\newcommand{\End}{\mathrm{End}}

\newcommand{\ggot}{\mathfrak{g}}

\newcommand{\Wc}{\mathcal{W}}

\newcommand{\Tr}{\mathrm{Tr}}

\newcommand{\Id}{\mathrm{Id}}

\newcommand{\la}{\langle}
\newcommand{\ra}{\rangle}

\renewcommand{\Im}{\mathrm{Im}\,}

\def\build#1_#2^#3{\mathrel{\mathop{\kern 0pt#1}\limits_{#2}^{#3}}}

\title{Integration formulas for Brownian motion on classical compact Lie  groups}

\author{Antoine Dahlqvist\thanks{ Statistical Laboratory, Centre for Mathematical Sciences, Wilberforce Road, Cambridge, CB3 0WA. E-mail address: ad814@maths.cam.ac.uk. }}

\begin{document}



\maketitle

\begin{abstract}Combinatorial formulas for the moments of the Brownian motion on classical compact Lie groups are obtained. These expressions are deformations   of formulas of B. Collins and P. \'Sniady for moments of the Haar measure and yield a proof of the First Fundamental Theorem of invariant theory and of classical Schur-Weyl dualities based on stochastic calculus.
\end{abstract}





\section{Introduction}
Let $G$ be a compact Lie group belonging to one of the three classical series: the unitary groups $\U(N),$ the orthogonal groups $\OO(N)$ and the compact symplectic groups $\Sp(N)$.  The aim of this article is  to give an  explicit expression to \begin{equation}
\int_G f(g)\mu(dg) \label{int}
\end{equation}
where $f$ is  a  polynomial function of the coefficients  of the matrix  and their adjoint and $\mu$ is a measure associated to the marginal  of a Brownian motion on $G$. 

\vspace{0,2 cm}

When $\mu $ is the Haar measure, such a study was started  by the physicist D. Weingarten (\cite{Weing}), a complete answer was achieved  later  on   by  B. Collins  in  ~\cite{ThCollins} for the unitary series and by  the same author together with  P. \'Sniady in  \cite{CS} for any of the above series.  For a more recent point of view, we refer  to  ~\cite{CM, ZJ}.  These results were obtained  using whether representation theory of the three series and the orthogonality of characters (\cite{ThCollins}),  or the first theorem of invariant theory or equivalently  Schur-Weyl dualities ~\cite{CS,CM,ZJ}.

\vspace{0,2 cm}

 When $\mu$ is a marginal of a Brownian motion,  this very question has been tackled by  T. Lévy  in ~\cite{ThierrySW},  following previous works \cite{IS,BianeMBL, Xu} focusing on polynomials of the form $\Tr(g^n),$ for $n\in \N.$   

\vspace{0,2 cm}

 The starting point of this article is to extend the results and the point of view of ~\cite{ThierrySW} in the following way: a first question answered here is to consider,  in the unitary case, not only polynomials of coefficients of matrices but also polynomials in the coefficients and their adjoints.  A second one is to get formulas that are both well suited for   large $N$ and large $t$ asymptotics.  Letting $t$ go to infinity, we recover formulas of  ~\cite{CS}, for the integration against the Haar measure. Another motivation for the present work is to answer the latter questions without using representation theory or the Theory of invariants.  Indeed, our proofs only rely on stochastic calculus and we moreover obtain a new proof of the First Fondamental Theorem of invariants theory (often abbreviated as FFT) and recover in this way the Schur-Weyl duality theorems for the classical compact Lie groups. Another feature of our work is that we express expectations  of functions measurable with respect to a Brownian motion on any  compact matrix Lie group,  in terms of expectations   with respect to the Brownian motion on unitary matrices of  the same size.

\vspace{0,2 cm}

The article is organized as follows. \textbf{Section 2} gives the definition and the scaling we are using for a Brownian motion on a classical compact Lie group. \textbf{Section 3} recalls the statements of the first fondamental theorem of invariants and of Schur-Weyl dualities and their equivalence. \textbf{Section 4} states and proves formulas for the Brownian motions. We prove here the formula for mixed moments of unitary matrices and their adjoint and recall the former results of T. Lévy. Then, we give an expression  for expectations of  tensors of Brownian motions on any classical group in terms of non mixed moments of the unitary Brownian motion. This is stated in Theorem \ref{Unitary representation}.  Using the expression of T. Lévy for non-mixed moments, we deduce from it a new combinatorial expression for expectation of entries of a Brownian motion for any of the three series, in terms of functions on the symmetric group. We obtain in this way the Theorem \ref{Deformed Weingarten}.  \textbf{Section 5} gives  two corollaries of the latter:  we obtain combinatorial formulas for integration against the Haar measure, that we compare with the one of  ~\cite{CS},  and a new proof of the first fondamental theorem of invariants. To show the convergence of our expressions as the time goes to infinity, we  have to show that expectations of non-mixed tensors of a Brownian motion vanish. This is simply proved using the factorization of a unitary Brownian motion as the product of a Brownian motion on $\SU(N)$ with an  Brownian motion on the circle. To conclude, in \textbf{section 6}, we prove in our framework, for the sake of completeness,  two already known lemmas  for the invariants of the symplectic groups.

\section{Brownian motion on classical compact Lie groups}

We discuss here about the definition of the Brownian motion we use all along the text. We have chosen here to define it as the solution of an SDE in the linear space of matrices.

\vspace{0,2 cm}

\noindent\textbf{Classical compact Lie groups:} For any positive integer $N$, we shall consider the following groups of matrices:  the orthogonal group $\OO(N)=\{O\in M_N(\R): {}^tOO=\Id \},$ the unitary group $\U(N)=\{U\in M_N(\C): U^*U=\Id\}$ and the special unitary group $\SU(N)= \{U:\U(N): \det(U)=1\},$ and as $N$ is even, the unitary symplectic group $\Sp(N)=\{S\in \U(N): {}^tSJS=J\}$, where $J$ is the matrix $J= \left(\begin{array}{cc}0 & I_{N/2}\\
-I_{N/2} & 0 \end{array}\right ).$ Their Lie algebras are denoted  by small gothic letters:  $\oN= \{X\in M_N(\R): {}^tX+X=0\},$ $\uN=\{X\in M_N(\C): X+X^*=0\}$ and $\su(N)=\{X\in \uN: \Tr(X)=0\}$, and $\usp=\{X\in \uN: {}^tXJ+JX=0\}$.  Let $G$ be  one of the above mentioned compact Lie groups  with Lie algebra $\ggot$.  All along this text, $\eps$ shall be $1$ if $G= \OO(N),$  $-1,$ if $G= \Sp(N)$ and $0$ otherwise.

\vspace{0,2 cm}

\noindent\textbf{Choice of a metric:} For any $x,y \in \End(\C^N)$, let us define  
$$\langle x,y \rangle=- \frac{N}{2^{|\eps|}}\Tr(xy).$$
The restriction of the real bilinear form $\la \cdot, \cdot\ra$  to $\ggot$ defines a $G$-invariant scalar product on the three series of Lie algebras we are considering.  The choice of normalization of the latter scalar product is motivated by the following remarks.  Let   $K$ be the Gaussian vector on $\ggot$ with covariance given by $\la\cdot,\cdot\ra.$   When $\ggot=\uN,$ the  random matrix $iK$ has the same law as a matrix from the Gaussian unitary ensemble\footnote{though, when $\ggot=\oN$ or $\usp$, $iK$ is not equal in law to a matrix of the  Gaussian orthogonal (resp. symplectic) ensemble.}.  On the other hand, as $N\to \infty,$  $\dim_\R(\ggot)\sim \frac{N^2}{2^{|\eps|}}$ and  according to the law of large numbers,  almost surely, $\frac{1}{\dim_\R(\ggot)}\|K\|^2\to1.$  Therefore, for any of the above classical series,    $\frac{1}{N}\Tr(KK^*)=\frac{2^{|\eps|}}{N^2}\|K\|^2\to 1.$  Furthermore, the Wigner theorem and its generalization to other invariant ensemble  (see \cite{Mehta}) show that the empirical measure of eigenvalues of $iK$ converges  weakly towards the semicircle law $   \sqrt{|4-x^2|_+} \frac{dx}{\pi}. $ 

\vspace{0,2 cm}

\noindent\textbf{Wrapping a Brownian motion of $\ggot$ on $G$:} Let us denote by $(K^\ggot_t)_{t\ge 0}$ the classical Brownian motion on the Euclidean space\footnote{we shall drop the symbol $\ggot$ when the context is clear.} $(\ggot,\la\cdot,\cdot\ra),$ that is the unique Gaussian process with covariance given for all $x,y\in \ggot$ and $ t,s\ge 0$, by $$\esp[\la x, K_t\ra \la y, K_s\ra ]= \min(t,s) \la x,y\ra.$$
The matrix  of quadratic variations $\la dK_t. dK_t \ra$ is equal to $C_\ggot dt $, where $C_\ggot$ is a deterministic matrix commuting with $G$.   Let  $S$ be  an element of $G$ and $(G_t)_{t\ge 0}$  be  the $\End(\C^N)$-valued stochastic process  solution to the following stochastic differential equation \begin{equation}
dG_t = G_tdK_t  +\frac{C_\ggot}{2}G_tdt, G_0=S.\label{BrownianWrapping}
\end{equation}

This Markovian process is called the \emph{Brownian motion} on $G$ issued from $S$.  When not precised, we shall assume that $(G_t)_{t\ge 0}$ is issued from $\Id.$ The rest of the paper being based on this definition, we shall  recall why it is justified. For any $x\in \ggot$, let $\Lc_x$ be the left-invariant first order differential operator defined by setting for any differentiable function $F\in C^\infty(\End(\C^N))$ and any $\phi\in \End(\C^N)$, 

$$\Lc_x F(\phi)=\left.\frac{d}{dt}\right|_{t=0} F(\phi\exp(t x) ) .$$
Let $\Uc(\ggot)$ be the real enveloping algebra of $\ggot$.  The application $\Lc$ extends to a morphism of associative algebra  between  ~$\Uc(\ggot)$ and the algebra of differential operators on $\End(\C^N)$ commuting with left translations by $G$.  For any orthonormal basis $(x_i)_{1\le i\le d}$ of $\ggot$ with respect to  the scalar product $\langle\cdot,\cdot\rangle$, let $c_\ggot=\sum_{i=1}^dx_i^2\in \Uc(\ggot)$. The element $c_\ggot$ does not depend on the orthonormal basis $(x_i)_{1\le i\le d}$ and is called the \emph{Casimir} element.  The image of the Casimir element $c_\ggot$ is denoted by $\Delta_G$: $$\Delta_G=\Lc_{c_\ggot}=\sum_{i=1}^d \Lc_{x_i}\circ \Lc_{x_i} .  $$
The image of $c_\ggot$ by the classical representation $\rho$ on $\C^N$ is $\sum_{i=1}^d \rho(x_i)^2= C_\ggot\in \End(\C^N)$. 
 The latter operators can be defined as well on $C^{\infty}(G),$ we denote them abusively by the same letters. The operator  $\Delta_G$ is then called  the \emph{Laplacian} associated to the metric  $\la\cdot,\cdot\ra$. 
\begin{lem} \label{Ito}For any $F\in C^\infty(\End(\C^N)),$  $(F(G_t))_{t\ge 0}$ is solution to the following SDE \begin{equation*}
d\left(F(G_t)\right)= \Lc_{dK_t}F(G_t)+\frac{1}{2} \Delta_\ggot (F)(G_t)dt.
\end{equation*}
\end{lem}
\begin{proof}For any $F\in C^\infty(\End(\C^N))$ and $\phi\in \End(\C^N),$ let  $dF_\phi\in \End(\C^N)^*, d^2_\phi F \in \left(\End(\C^N)^{\ts 2}\right)^* $ be the one form and two  form associated to the first and second derivatives of $F$. For any  $x\in \ggot, \phi\in \End(\C^N)$,
 \begin{equation*}
\Lc_{x}\circ \Lc_x (F)(\phi) =d_\phi  F(\phi x^2) +d_\phi^2F(\phi x\ts \phi x).\tag{*}\label{second derivative}\end{equation*}
Let $(x_i)_{i=1}^d$ be an orthogonal basis of $(\ggot,\la\cdot,\cdot\ra),$ the Itô formula combined with the latter equation implies that 
\begin{align*}
d\left(F(G_t)\right)&=d_{G_t} F (G_t dK_t) +\frac{1}{2}  d_{G_t} F(G_t C_\ggot)dt +\frac{1}{2}\la d^2_{G_t} F ( G_t dK_t\ts G_t dK_t ) \ra\\
&= \Lc_{dK_t}F(G_t)+ \frac{1}{2}  d_{G_t} F(G_t C_\ggot)dt +\frac{1}{2} \sum_{i=1}^d d^2_{G_t} F ( G_t x_i\ts G_t x_i ) dt 
\end{align*}Combined with (\ref{second derivative}), the latter equation yields the result.\end{proof}
It can be proved  that  the process $(G_t)_{t\ge 0}$ satisfies the following properties (more details can be found in \cite{Liao}).
\begin{lem}\label{PropMB}i) Almost surely,  for any $t\ge 0,$ $G_t\in G.$  \\
ii) For all $T\ge 0,$  $(G_T^{-1}G_{T+t})_{t\ge 0}$ is independent of $\sigma(G_s,s\le T)$ and has the same law as $(gG_tg^{-1})_{t\ge 0}$ for any $g\in G$.\\
iii) Almost surely, the mapping $t\in\R_+\mapsto G_t\in G$ is continuous.\\
iv) The generator of the Markov process $(G_t)_{t\ge 0}$ on $G$ is $\Delta_G.$
\end{lem}
\begin{proof} We only sketch the first point and leave the others to the Reader (see \cite{Liao} for proofs).  For any of the above groups, there is a function  $F\in C^\infty(\End(\C^N),\R^n)$ such that $G= F^{-1}(\{1\}).$  For any $x\in \ggot, $ $\Lc_x F =\Delta_\ggot F=0$. Hence, lemma \ref{Ito} yields that almost surely for any $t\ge 0,$ $F(G_t)=1.$ 
\end{proof}

For $G=\U(N),$ the Brownian motion satisfies the following splitting property.
\begin{lem}\label{Facto U SU Cercle} Let $(U_t)_{t\ge 0}$ be a Brownian motion on $\U(N)$.  Then $(U_t)_{t\ge 0}$ has the same law as $(e^{\frac{iB_t}{N}}S_t)_{ t \ge 0}$, where $(B_t)_{t\ge 0}$ is a standard real Brownian motion and $(S_t)_{t\ge 0}$ an independent Brownian motion on $\SU(N)$. 
\end{lem}
\begin{proof} Let us set for any $t\ge 0,$ $\tilde{K}_t= K^{\su(N)}_t+\frac{1}{N} i B_t,$ where $(B_t)_{t \ge0}$ is a standard real Brownian motion independent from $( K^{\su(N)}_t)_{t\ge 0}.$  On the one  hand,   $(\tilde{K}_t)_{t\ge 0}$ has the same law as $(K^{\uN}_t)_{t\ge 0}.$ On the other hand,  since $(S_t)_{t\ge 0 }$ satisfies (\ref{BrownianWrapping}) with driving noise  $( K^{\su(N)}_t)_{t\ge 0},$  so does $(e^{\frac{iB_t}{N}}S_t)_{ t \ge 0}$ but with driving noise $(\tilde{K}_t)_{t\ge0}.$
\end{proof}

Our aim is to get  a formula for $\esp[P(G_t)],$ for any homogeneous polynomial $P$ in coefficients  of matrices and their adjoint, that is well suited to let $t\to \infty$ and to get another proof of the formulas of ~\cite{CS} for $\esp[P(H)],$ where $H$  is the Haar measure on $G$.  We are going to see that for any series, these expectations can be expressed in terms of expectations of coefficients of unitary Brownian motions (and not their adjoint).  For a fix degree of homogeneity, we shall  gather these expectations in a single one, namely $\esp[ G_t^{\ts n}\ts \overline{G_t}^{\ts m}],$ that will be  expressed in terms of its Schur-Weyl dual. The latter expression will imply the FFT theorem and  therefore the Schur-Weyl duality.  Before proceeding,  we recall  the statement of FFT and its relation to Schur-Weyl duality.

\section{First Fondamental Theorem of invariants and Schur-Weyl duality\label{SectionFFT}}

For any vector space $W$ on which $G$ acts, we consider the set $W^G$ of points  fixed by $G.$ Let us write $V$ for the canonical representation $\C^N.$ The fondamental theorem of invariant theory gives a spanning family for the vector space $W^G$, when $W$ is a space of tensors of the natural representation $V$ or its dual $V^*$.  We refer to ~\cite{GW,Proc} for a nice exposition and a proof of this theorem relying on algebraic geometry.  Let us recall its statement. We will  restrict  to the three series of groups $\OO(N),\U(N)$ and $\Sp(N),$ discarding here the special unitary groups $\SU(N)$.   For any $n,m\in \N,$ we shall consider the representation $(V^{\ts n}\ts {V^*}^{\ts m},\rho_{n,m})$ of $G$   defined by  $$\rho_{n,m}(g). v_1\ts \ldots \ts v_n \ts \vp_1\ts \ldots \ts \vp_m=(g. v_1)\ts \ldots \ts (g.v_n) \ts \vp_1(g^{-1}\cdot)\ts \ldots \ts \vp_m(g^{-1}\cdot),$$
for all $g\in G, v_1,\ldots,v_n\in V, \vp_1,\ldots, \vp_m\in V^*.$  Let us highlight that when $G=\U(N),$ $V^*$ is isomorphic as a representation to $\overline{V},$ that is, the real vector space $V,$ endowed with the actions $\lambda$ and $\rho$  of $\C$ and $\U(N),$ defined  respectively  by $\lambda(z,v)=\overline{z}v $ and $\rho(U, v)= \overline{U}v, $ for any $z\in \C,$ $U\in \U(N)$ and $v\in V.$  When $G$ is respectively  $\OO(N),$  $\Sp(N)$ and $\U(N),$  we denote by $B$   the canonical $G$-invariant $\C$-bilinear  form on  $V\oplus V^*$, that is   skew-symmetric as $G=\Sp(N)$ and symmetric when $G$ is $\OO(N)$ or $\U(N).$  Let  $\Mc(n)$ denote the set of \emph{matchings} of $\{1,\ldots,2n\},$ that is, partitions of $\{1,\ldots, 2n\}$ whose blocks are all of size $2$.

 \begin{theo}[FFT] \label{FFT} The vector space  $\left(V^{\ts n}\ts {V^*}^{\ts m}\right)^G$ is trivial if  $n+m$ is odd, or $G=\U(N)$ and $n\not=m$. Moreover, as $n+m$ is even, it is spanned by the linear maps 
 $$\begin{aligned} {V^*}^{\ts n}\ts V^{\ts m}&\longrightarrow \C \\(\phi_1,\ldots, \phi_{n+m})&\longmapsto \prod_{\{a,b\}\in\pi} B(\phi_a,\phi_b), \hspace{1 cm}\pi\in \Mc(\frac{n+m}{2}).\end{aligned}$$
  \end{theo}
  \begin{rema} When $n=dN$ and $\U(N)$ is replaced by $\SU(N),$ or when $N$ is even, $n=\frac{N}{2}$ and $\OO(N)$ is replaced by $\SO(N),$ the Theorem \ref{SW} does not hold anymore, invariants associated to the determinant have to be added (see ~\cite{Proc,Brauer37}).  
\end{rema}

The first part of the Theorem can be proved with the following remarks.  Let us consider $w\in  \left(V^{\ts n}\ts {V^*}^{\ts m}\right)^G$.When  $n+m$ is odd, then   $-\Id\in G$ and $w=-w$. When $G=\U(N)$ and $n\not=m,$  $z\Id\in \U(N)$ for any $z\in \U(1),$ and $w= z^{n-m}w.$

\vspace{0,2 cm}

In the latter theorem, when the invariant space is non trivial, the vector space $V^{\ts n}\ts {V^*}^{\ts m}$ can be identified with a space of endomorphisms, and  invariant vectors with endomorphisms commuting with $G$.      Indeed, when $G$ is $\OO(N)$ or $\Sp(N),$  the map $$\begin{aligned}\theta: V&\longrightarrow V^*\\v&\longmapsto B(v,\cdot) \end{aligned}$$ is a linear isomorphism of representation, and if $n+m$ is even, $V^{\ts n}\ts {V^*}^{\ts m}$ is isomorphic to $\End(V^{\ts \frac{n+m}{2}})$. For any vector space $W$ endowed with an action of $G$, let us denote the commutant of $G$ in $\End(W)$ by
$$ \End_G(W)=\{f\in \End(W): g. f(w) = f(g.w), \text{ for any }w\in W,g\in G\}.$$
Using the application $\theta: V\to V^*$, the first theorem of invariant theory yields a spanning family of the vector space $\End_G(V^{\ts \frac{n+m}{2}})$. These statements are known as \emph{Schur-Weyl dualities}.  The latter space is furthermore an algebra that we shall describe in the next paragraph.  For the symplectic groups, we need there to fix an appropriate sign for each invariant. A natural one is given for invariants  in $\left(V^{\ts 2n}\right)^{\Sp(N)}.$  Let us therefore spell out  first the statement of Theorem \ref{FFT} for $\left(V^{\ts 2n}\right)^{\OO(N)}, \left(V^{\ts n}\ts \overline{V}^{\ts n}\right)^{\U(N)}$ and $\left(V^{\ts 2n}\right)^{\Sp(N)}$. \vspace{0,2 cm}

\noindent\emph{Indexing invariants by $\Sy_{2n}/\Hc_n$:} Let $(e_i)_{1\le i\le N}$ be the canonical basis of $V$.     Two combinatorial sets allow to describe the latter invariants: the set $\Mc(n)$ of matchings of $\{1,\ldots,2n\}$ and its subset  $\Mc_+(n)$, consisting of elements matching $\{1,\ldots,n\}$ with $\{n+1,\ldots,2n\}.$  Let us set $\pi_0= \{\{k,2n+1-k\},k\in\{1,\ldots,n\}\}\in \Mc_+(n).$ The following element is a  $G$-invariant
 $$w_{\pi_0}= \sum_{1\le i_1,\ldots,i_{2n}\le N}  \prod_{k=1}^n B(e_{k},e_{2n+1-k}) e_{i_1}\ts \ldots \ts e_{i_{2n}}.$$
Theorem \ref{FFT} states that all invariants are obtained by permutations of  the tensors.  Let us highlight a choice of equivalence class representatives.   Both $\Mc(n)$ and $V^{\ts 2n}$ (resp. $\Mc_+(n)$ and $V^{\ts n }\ts \overline{V}^{\ts n}$) are endowed with a left action of the permutation group $\Sy_{2n}$ (resp. $\Sy_n\times \Sy_n$).  Let us recall that the second one is  defined  by setting for any  $\sigma \in\Sy_{2n},$ and $v\in V^{\ts 2n},$ 
\begin{equation}
\sigma. v =  v_{\sigma^{-1}(1)}\ts \ldots\otimes v_{\sigma^{-1}(2n)}.\label{actionGpSym}
\end{equation}
 Let  $\Hc_n$ (resp. $\Dc_n$) be the stabiliser of $\pi_0$  in $\Sy_{2n}$ (resp. $\Sy_n\times \Sy_n$), so that $\Mc(n)= \Sy_{2n}/\Hc_n$ (resp. $\Mc_+(n)=\Sy_n\times\Sy_n /\Dc_n$). It is isomorphic to the  \textit{hyperoctahedral group} (resp. $\Sy_n$). What is more,  for any $h\in \Hc_n,$  $h .w_{\pi_0}=w_{\pi_0},$ when $G$ is $\U(N)$ or $\OO(N)$, whereas $h .w_{\pi_0}= \eps(h)w_{\pi_0}, $ when $G=\Sp(N),$  where $\eps:\Sy_{2n}\to\{-1,1\}  $ denotes the signature morphism.  This consideration shows that the following elements are well defined: for any $\sigma\in \Sy_{2n},$ if $G=\OO(N)$, or   $G=\U(N)$ and $\sigma \in \Sy_n\times \Sy_n$,  let us set 
\begin{equation}
w_{\sigma.\pi_0}= \sigma .w_{\pi_0},\label{Spanning Family of Invariant}
\end{equation}
and if  $G=\Sp(N),$ set 
\begin{equation}
w_{\sigma.\pi_0}=\eps(\sigma)\sigma .w_{\pi_0}.\label{signSymp}
\end{equation}

\noindent\emph{Brauer algebras:} Let   $\Bc_n$  be the vector space  of functions on   $\Mc(n)$. We write $\rho:\Sy_n\to \End(V^{\ts n})$ the representation of the symmetric group given by permutation of tensors, as defined in (\ref{actionGpSym}),  and  $\rho_{+},\rho_{-}: \Bc_n\to \End(V^{\ts n})$      the linear maps such that 
$$\rho_{\eps}(\pi)=\Id^{\ts n}\ts {\theta}^{\ts n} (w_\pi),$$
where we have identified $V^{\ts n}\ts {V^{*}}^{\ts n}$ with $\End(V^{\ts n})$ and the right-hand-side is resp. given by (\ref{Spanning Family of Invariant}) and (\ref{signSymp}), when  the symbol $\eps$ is resp. $+$ and $-$.   For $n,m\in \N^*,$ let us define $s:\{1,\ldots, 2(n+m)\}\to \{-1,1\},$ with $s(k)= -1$ if $n<k\le n+2m$ and $1$ otherwise. We  shall also consider the vector space  $\Bc_{n,m}$ with basis  indexed by matchings $\pi$ of $\{1,\ldots, 2(n+m)\}$  with $s(a)\not=s(b)$ for any block $\{a,b\}$ of $\pi$. The  Schur-Weyl duality theorem is an easy consequence of the theorem  \ref{FFT}, we leave this proof as an exercice for the Reader (see for example Section 10 of ~\cite{GW}).

\begin{figure}[!h] 
  \centering
 \includegraphics[height=2in]{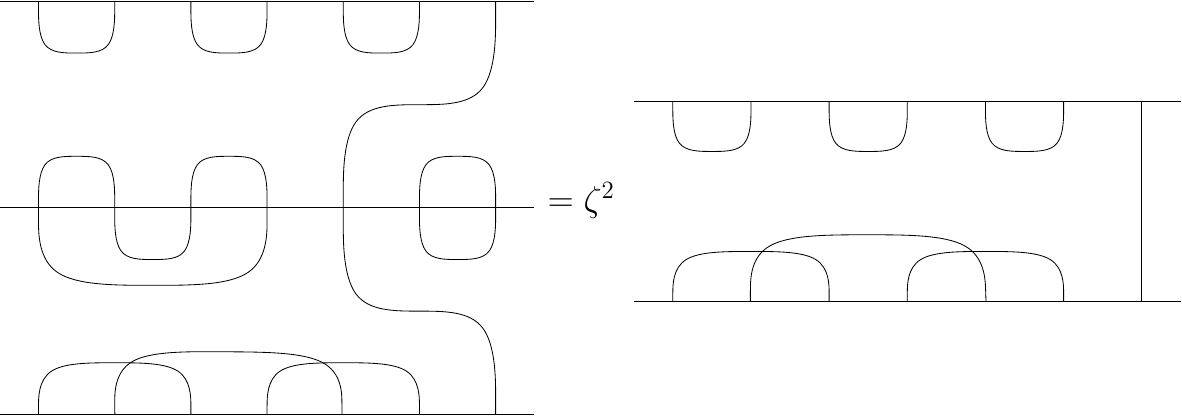}
  \caption{Multiplication  $\pi.\eta$,   with $\pi=\{\{1,2\},\{3 ,4\},\{ 5,6 \},\{7 ,10\},\{ 8,9\},\{ 11,12\},\{13 ,14\}\}$ and $\eta=\{\{1, 4\},\{2 , 3\},\{ 5,8 \},\{6 ,7\},\{ 9,11\},\{ 10,13\},\{12 ,14\}\}$.}\label{multB}
\end{figure}

\begin{figure}[!h]
\centering
\includegraphics[height=1in]{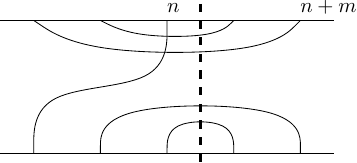}
\caption{An element of  $\Bc_{n,m}(\zeta)$.}\label{BrauerMur}
\end{figure}

\begin{theo}[Schur-Weyl Duality] \label{SW}The action of $G$ on $V^{\ts n}$   has the following commutant:
$$\End_G(V^{\ts n})= \left\{\begin{array}{ccc} \rho(\C[\Sy_n]), &\text{ if }G=\U(N),\\\rho_+(\Bc_n), &\text{ if }G=\OO(N),\\\rho_-(\Bc_n), &\text{ if }G=\Sp(N).\end{array}\right.
$$
Furthermore, for $n,m\in \N^*,$ $$\End_{\U(N)}(V^{\ts n}\ts \overline{V}^{\ts m})= \rho_+(\Bc_{n,m}).$$
\end{theo}
A nice property of the applications $\rho,\rho_+$ and $\rho_-$ is that there exists  structures of algebra on $\Bc_n$, such that they are algebra morphism.  Let us fix a complex number $\zeta\in \C.$ Let us represent matchings in  $\Mc(n)$  as simple  curves in the plane with endpoints     $(-\frac{n-1}{2},\frac{1}{2}),(-\frac{n-1}{2}+1,\frac{1}{2}),\ldots, (\frac{n-1}{2},\frac{1}{2})$  and $(\frac{n-1}{2},-\frac{1}{2}), (\frac{n-1}{2}-1,-\frac{1}{2}),\ldots, (-\frac{n-1}{2},-\frac{1}{2}),$ by labelling endpoints in this order from $1$ to $2n$. Given two matchings $\pi$ and $\nu\in \Mc(n)$, the concatenation of the curves of $\pi$ translated by $(0,1)$ with the curves of $\nu$ yields curves with endpoints  $(-\frac{n-1}{2},\frac{3}{2}),(-\frac{n-1}{2}+1,\frac{3}{2})\ldots, (\frac{n-1}{2},\frac{3}{2}),$ $(\frac{n-1}{2},-\frac{1}{2}),\ldots,(-\frac{n-1}{2},-\frac{1}{2}),$ which  represents  a  matching $\pi\circ \nu\in\Mc(n),$ and $b(\pi,\nu)$  loops. We set\footnote{see figure \ref{multB} for an example.} $$\pi.\nu=\zeta^{b(\pi,\nu)}\pi\circ \nu.$$ 
This is a simple exercice to check that the extension of this operation by linearity defines a multiplication on $\Bc_n.$   This algebra is denoted by $\Bc_n(\zeta)$ and  called the \emph{Brauer algebra} (see \cite{Brauer37} or \cite{GW}[Section 10.1]).     
For any $n,m\in\N^*,$ a matching of $\Mc(2(n+m))$ belongs to $\Bc_{n,m}$   if and only if it can be represented  as a collection of simple curves that cross exactly once  the lines $\R\times \{0\}$ or $\{\frac{n-m}{2}\}\times\R $ (see figure \ref{BrauerMur}).  Let us remark that for all $x,y \in \Bc_{n,m},$   their product in $\Bc_{n+m}(\zeta)$ satisfies $x.y\in\Bc_{n,m}.$ Therefore, this operation defines an algebra structure, we then denote this algebra by $\Bc_{n,m}(\zeta).$ It is called the \emph{walled Brauer algebra} (it has been originally considered in  \cite{WBTuraev,Koike} , see also \cite{WBrauerHalverson}). Note that the element $\pi_0$ is the identity of these algebras, we shall denote it alternatively by $1$.

\begin{lem} \label{algebra morphisms}For any $n,m\in\N^*,$ the applications $\rho_+: \Bc_n(N)\to \End_{\OO(N)}(V^{\ts n}),$   $\rho_-: \Bc_n(-N)\to \End_{\Sp(N)}(V^{\ts n})$ and $\rho_+:\Bc_{n,m}(N)\to\End_{\U(N)}(V^{\ts n}\ts \overline{V}^{\ts m}) $ are  onto algebra morphisms.
\end{lem}
We shall leave the proof of this fact as an easy exercice for the unitary and orthogonal case and give one in the section  \ref{Proof Sp} for the symplectic case,  as we have not found a direct proof in the literature.
For any $1\le i,j\le \dim(V),$ let $E_{i,j}$ denote the endomorphism of $\End(V),$ such that  $E_{i,j}(e_k)=\delta_{j,k}e_i,$ for any $1\le k\le \dim(V).$ Then, for any $\pi\in \Mc(n),$  
\begin{equation}
\rho_+(\pi)= \sum_{1\le i_1,\ldots,i_{2n}\le N} \left(\prod_{\{a,b\}\in\pi}\delta_{i_a,i_b}\right) E_{i_1,i_{2n}}\ts E_{i_2,i_{2n-1}}\ldots \ts E_{i_n,i_{n+1}}.\label{RepO}
\end{equation}
whereas if $\pi= \sigma.\pi_0,$ with $\sigma\in \Sy_{2n},$
\begin{equation}
\rho_-(\pi)= \sum_{1\le i_1,\ldots,i_{2n}\le 2N} \eps(\sigma) \left(\prod_{\{a,b\}\in\pi_0}J_{i_{a},i_b}\right) E_{i_{\sigma^{-1}(1)},i_{\sigma^{-1}(2n)}}\ts\ldots \ts E_{i_{\sigma^{-1}(n)},i_{\sigma^{-1}(n+1)}}{J^{-1}}^{\ts n}.\label{RepSp}
\end{equation}

To close this section, let us highlight two kinds of matchings that will play thereafter a prominent role. For any  integers $1\le a<b\le n$, we set

$$\la a\,b \ra=\{ \{a ,b\}, \{2n+1-a, 2n+1-b\} \}\cup\{\{k,2n+1-k\}, 1\le k\le n, k\not\in\{a,b \}  \}$$   and  the usual transposition is identified as 
$$(a\, b)=\{ \{a,2n+1-b\},\{ b, 2n+1-a\}\}\cup\{\{k,2n+1-k\}, 1\le k\le n, k\not\in\{a,b \}  \}.$$
\begin{figure}[!h]
\centering
$\begin{array}{ccc}
 \includegraphics[height=1in]{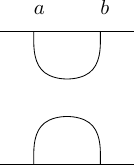} &\hspace{2 cm} & \includegraphics[height=1in]{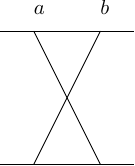}   \end{array}$
\caption{The two elements $\la a\,b\ra$ et $(a\,b)$ of $\Bc_n(\zeta)$.}\label{TranspoContrac}

\end{figure}

\noindent When $n=2,$ they  have the following linear representations  on $V^{\ts 2},$
\begin{equation}
\rho_+(   \la1\, 2 \ra )=  \sum_{1\le i,j\le N} E_{i,j}\ts E_{i,j}, \hspace{0,5 cm}\text{ and}\hspace{1 cm}\rho_+(   (1\, 2 ) ) = \sum_{1\le i,j\le N} E_{i,j}\ts E_{j,i}, \label{RepContracO}
\end{equation}
whereas 
\begin{equation}
\rho_-(   \la1\, 2 \ra )=  -\sum_{1\le a,b,c,d\le N}J_{a,c}J_{b,d} E_{a,b}\ts E_{c,d}, \hspace{0,5 cm}\text{ and}\hspace{1 cm}\rho_-(   (1\, 2 ) ) = -\sum_{1\le i,j\le N} E_{i,j}\ts E_{j,i}.  \label{RepContracSp}
\end{equation}
For $n\ge 2$ and any $1\le a<b\le n$, $\rho_\eps(\la a\,b\ra)$ acts like $\rho_\eps(\la 1\,2\ra)$  on the $a$-th and $b$-th factors of $V^{\ts n}$ and trivially elsewhere. Let us lastly remark that the family $\{(a\, b),\la a\,b\ra: 1\le a<b\le n\}$ generates $\Bc_n(\zeta)$ as an algebra. 




\section{Integration formulas for the Brownian motion}

We shall here recall the result of ~\cite{ThierrySW} about integration against a Brownian motion and give a slight extension of it for the unitary case. Starting from these formulas, we show that any moment (possibly mixed) for any series  can be expressed in terms of non-mixed moments for the unitary groups. The  constant  (in time)  part of this expression corresponds  to the integration against the Haar measure. We then check that these expressions do match the one of ~\cite{CS}.
\begin{prop}[\cite{ThierrySW}, Prop. 2.2, 2.6 and 2.8] For any $t\ge 0,$    $\esp[G_t^{\ts n}]=e^{t \Delta_{G}(n)},$ where
$$  - \Delta_{G}(n)= 
\left\{\begin{array}{ll}\rho\left(\frac{n}{2}+\frac{1}{N}\sum_{1\le i<j\le n}(i\, j) \right),&\text{ if } G=\U(N),\\
\rho_{\eps}\left(\frac{n}{2}(1-\frac{1}{\eps N})+\frac{1}{\eps N}\sum_{1\le i<j\le n}((i\, j)  -\la i\,j\ra    )\right), &\text{ if } G=\OO(N) \text{ or }\Sp(N).
\end{array}\right. $$\label{BrownianSW}
\end{prop}
Our notations being  different from the one of \cite{ThierrySW}, we recall for completeness a proof for the symplectic case in the appendix (section \ref{Proof Sp}).
\begin{lem}  For $G=\U(N)$ and  any $n,m\in \N^*,$ $\esp[G_t^{\ts n}\ts \overline{G}_t^{\ts m}]= e^{t\Delta_{U}(n,m)},$ with
$$-\Delta_U(n,m)= \rho_+\left(\frac{n}{2}+\frac{1}{N}\sum_{1\le i<j\le n+m }(\un_{\{i,j\le n \text{ or }i,j>n\}}(i\, j)  -\un_{\{i\le n<j\}}\la i\,j\ra    )\right).$$
\end{lem}
\begin{proof} Let us recall that for any finite dimensional representation $(V_f,f)$,  setting for all $x\in \ggot,$ $f(x)=\Lc_x(f)(\Id)$, for any  $a,b\in \ggot,$   $\Lc_{a}\circ\Lc_bf(\Id)= f(a)f(b)\in\End(V_f). $ For example, 
  if $(x_i)_{1\le i\le \dim(\ggot)}$ is an orthogonal basis of $(\ggot,\la\cdot,\cdot\ra),$
$\Delta_G(f)(\Id)=\sum_{i=1}^d f(x_i)^2=f(c_\ggot).$ Therefore, for any $t\ge 0$, 
$$\frac{d}{dt}\esp[f(G_t)]=\frac{1}{2} \esp[\Delta_G(f)(G_t)]=\frac{1}{2}\esp[f(G_t)]\Delta_G(f)\in \End(V_f),$$
and   

$$\esp[ f(G_t)]=\exp(\frac{t}{2}f(c_\ggot)). $$
We need here to consider the representation $(V^{n}\ts {V^*}^{\ts m},\rho_{n,m}),$ identifying $V^*$ with $\overline{V}$. Let us  choose  the orthonormal basis  $\{ \frac{1}{\sqrt{2N}}(E_{k,l}-E_{l,k}),\frac{i}{\sqrt{2N}}(E_{k,l}+E_{l,k}): 1\le k<l\le N\}\cup\{\frac{i}{\sqrt{N}} E_{k,k}:1\le k\le N\}$ of $(\uN,\langle\cdot,\cdot \rangle)$ and compute the Casimir $c_{\uN}$as an element of the real  algebra $\Uc(\glN(\C))$:
\[2N c_\uN=\sum_{1\le k,l\le N}  (1+ i\ts i ) E_{k,l}\ts E_{k,l}+ ( i\ts i -1) E_{k,l}\ts E_{l,k}\in\Uc(\glN(\C)). \]
Considering  $f=\rho_{n,m}$ as a representation of the real  enveloping algebra $\Uc(\glN(\C))$  implies that 
\begin{align*}
2N \rho_{n,m}(c_\uN)&= (n+m)N\Id_{V^{\ts n+m}}+2\sum_{1\le a\le n<b\le n+m} \langle a \hspace{0,2cm} b\rangle-2\sum_{1\le a<b\le n, \text{ or }n<a<b\le n+m }( a\hspace{0,2cm} b).
\end{align*}
\end{proof}
Note that at first sight, these combinatorial expressions do not shed any light  about the integration against the Haar measure. We give in what follows  a reformulation of the latter that covers up this point.  Our formula uses the  following notations.

\vspace{0,2 cm}

\noindent\emph{Partial  matchings:}  We shall consider the sets   $\Pc\Mc(n)$ of partitions of  $\{1,\ldots,n\}$, with blocks of size $1$ and $2$.    For $n,m\in\N^*$, we denote by $\Pc\Mc(n,m)$ the subset of $\Pc\Mc(n+m)$ of elements without  $2$-blocks connecting $\{1,\ldots,n\}$ with itself. For any $\Pi$ in $ \Pc\Mc(n),$ we denote by $|\Pi|$ its number of $2$-blocks. 

\vspace{0,2 cm}

\noindent\emph{Matched  tensors:} For $a,b,c\in \{1,\ldots, n\},$ with $b<c,$ let us define two linear maps from respectively $\End(V)$ and $\End(V)^{\ts 2}$ towards $\End(V^{\ts n}),$ by setting for any $A,B,C\in\End(V),$
$$(A)_a= \Id^{\ts a-1}\ts A\ts \Id^{\ts n-a}$$
and 
$$(B\ts C)_{b,c}= \Id^{\ts b-1}\ts B\ts \Id^{\ts c-b-1}\ts C\ts \Id^{n-c}.$$
Given a  partial matching $\pi\in \Mc(n)$,  let us define a  map from $ \End(V)\times\End(V)^{\ts 2}$ to $\End(V^{\ts n})$, setting for any $A\in\End(V)$ and  $T\in \End(V)^{\ts 2},$ 
$$(A,T)^{\ts \pi}=\prod_{\{i\}\in\pi } (A)_i  . \prod_{\{a,b\}\in \pi: a\not=b} \left(T\right)_{a,b}\in \End(V^{\ts n}). $$
Given a partial matching $\pi\in \Pc\Mc(n,m),$ $A,B\in \End(V)$ and $T\in \End(V)^{\ts 2}$, we set 
$$\left((A, B), T\right)^{\ts \pi }=\prod_{\{i\}\in\pi :i\le n} (A)_i.\prod_{\{i\}\in\pi: i>n } (B)_i   . \prod_{\{a,b\}\in \pi: a\not=b} \left(T\right)_{a,b}\in \End(V^{\ts n}). $$

\noindent Let $(U_t)_{t\ge 0}$ and $(V_t)_{t\ge 0}$ be two independent $\U(N)$-Brownian  motions. 
\begin{theo}For any $t\ge 0,$  when $G$ is $\OO(N)$ or $\Sp(N),$
$$\begin{aligned}\esp[G_t^{\ts n}]=\sum_{\pi\in\Pc\Mc(n)} \esp\left[ \left(e^{\frac{\eps t}{2N}}U_t, \int_{0}^te^{\frac{ s\eps}{N}}U_s^{\ts 2}ds \right)^{\ts \pi}  \right] \prod_{ \{a,b\}\in \pi }\frac{  \rho_{\eps}(\la a\, b\ra )}{ \eps N}\\
\end{aligned}.$$
When $G=\U(N),$
$$\esp[U_t^{\ts n}\ts \overline{U}^{\ts m}_t]=\sum_{\pi\in \Pc\Mc(n,m)}  \esp\left[ \left((U_t, V_t),\int_{0}^tU_s\ts V_sds \right)^{\ts \pi}  \right]\prod_{ \{a,b\}\in \pi }  \frac{ \rho_{+}(\la a\, b\ra )}{N}.$$\label{Unitary representation}
 \end{theo}
 As we will see in the next section, these formulas easily yield an expression for the integration against the Haar measure. As $t\to\infty,$ the only terms remaining are the one indexed by partial  matchings with no block of size $1$, except in the orthogonal case. 
  \begin{proof} We shall consider the case $G\in \{\OO(N),\Sp(N)\},$ the third one being very similar.  We will argue that  the function $R$ on the right-hand-side  satisfies $\frac{d}{dt }R_t= R_t\Delta_G(n)\in \End(V^{\ts n})$ and $R_0=\Id.$ According to Proposition  \ref{BrownianSW}, these equations are true for the left-hand-side and this therefore yields the announced equality. First note that when $t=0,$ the matched tensors indexed by $\pi\in \Pc\Mc(n)$ appearing in the right-hand-side vanish   if $|\pi|\not=0$ and equal $\Id$  otherwise. Hence, the initial condition holds true.   For any partial matching $\pi\in\Pc\Mc(n)$,  $\{a,b\}\in \pi$  and $\{c\},\{d\}\in\pi,$ with $c\not=d,$ we shall denote respectively by $\hat{\pi}_{a,b}\in\Pc\Mc(n)$ and $\check{\pi}^{c,d},$the partial matchings with    $\{a,b\}$, resp. $\{c\},\{d\} ,$   replaced by $\{a\},\{b\}$,  resp. $\{c,d\}$.  On the one hand, Itô's formula implies   
 \begin{equation}
\begin{aligned}\frac{d}{dt}\esp&\left[ \left(e^{\frac{\eps t}{2N}}U_t, \int_{0}^te^{\frac{ s\eps}{N}}U_s^{\ts 2}ds \right)^{\ts \pi}\right]= \sum_{\{a,b\}\in\pi}\esp\left[\left(e^{\frac{\eps t}{2N}}U_t, \int_{0}^te^{\frac{ s\eps}{N}}U_s^{\ts 2}ds \right)^{\ts \hat{\pi}_{a,b}}\right]\\&\hspace{2 cm}-\esp\left[ \left(e^{\frac{\eps t}{2N}}U_t, \int_{0}^te^{\frac{ s\eps}{N}}U_s^{\ts 2}ds \right)^{\ts \pi}\right]\left(\frac{n}{2}\left(1-\frac{\eps}{N}\right) +\frac{1}{N}\sum_{\{a\},\{b\}\in\pi:\,  a<b } (a\,b)\right).\end{aligned}\label{Diff Matched Tensors}
\end{equation}
On the other hand, for any integer $1\le h\le \frac{n}{2}$ and any $t\ge 0,$  $\prod_{k= 0}^{h-1} \la n-2k\,\, n-2k-1 \ra \esp[G_t^{\ts n}]= \prod_{k= 0}^{h-1} \la n-2k\,\, n-2k-1 \ra \esp[G_t^{\ts n- 2h}\ts \Id^{\ts 2h}]$.    Differentiating this equality yields 
\begin{equation}
\prod_{k= 0}^{h-1} \la n-2k\,\, n-2k-1 \ra \Delta_{G}(n)= \prod_{k= 0}^{h-1} \la n-2k\,\, n-2k-1 \ra \Delta_{G}(n-2h)\ts \Id^{\ts 2h},\label{Projective family Casimir}
\end{equation}
 whereas conjugating by elements of $\Sy_{n}$ and using the expressions of Proposition \ref{BrownianSW} implies that for any partial matching $\pi \in \Pc\Mc(n),$
\begin{equation}
\begin{aligned}\prod_{ \{a,b\}\in \pi }\frac{  \rho_{\eps}(\la a\, b\ra )}{ \eps N}& \Delta_{G}(n)=\sum_{\{c\},\{d\}\in \pi:\,c<d}\prod_{\{a,b\}\in \check{\pi}^{c,d}}\frac{\rho_\eps(\la a\,b\ra)}{\eps N} \\& -\prod_{ \{a,b\}\in \pi }\frac{  \rho_{\eps}(\la a\, b\ra )}{ \eps N} \rho_\eps\left(\left(\frac{n}{2}-|\pi|\right)\left(1-\frac{\eps}{N}\right) +\frac{1}{N}\sum_{\{c\},\{d\}\in\pi:\,  c<d } (c\,d)\right).\end{aligned}\label{Explicit Projective family Casimir}
\end{equation}
Let us denote   for any partial matching $\pi\in \Pc\Mc(n),$  
\begin{equation}
-\Delta(\pi)=\left(\frac{n}{2}-|\pi|\right)\left(1-\frac{\eps}{N}\right) +\frac{1}{N}\sum_{\{c\},\{d\}\in\pi:\,  c<d } (c\,d). \label{NotLaplacienAppar}
\end{equation}
 For any $\pi\in \Pc\Mc(n),$ $\Delta(\pi)$ commutes with $\prod_{\{a,b\}\in\pi}\rho_{\eps}(\la a\, b\ra )$, so that using  (\ref{Diff Matched Tensors})  and  (\ref{Explicit Projective family Casimir}) implies 
$$\begin{aligned}\frac{d}{dt}R_t&= \sum_{\pi \in\Pc\Mc(n)}\sum_{\{c,d\}\in\pi} \esp\left[ \left(e^{\frac{\eps t}{2N}}U_t, \int_{0}^te^{\frac{ s\eps}{N}}U_s^{\ts 2}ds \right)^{\ts \hat{\pi}_{c,d}}\right]\prod_{ \{a,b\}\in \pi }\frac{  \rho_{\eps}(\la a\, b\ra )}{ \eps N}\\&-\sum_{\pi \in\Pc\Mc(n)} \esp\left[ \left(e^{\frac{\eps t}{2N}}U_t, \int_{0}^te^{\frac{ s\eps}{N}}U_s^{\ts 2}ds \right)^{\ts \pi}\right]\prod_{ \{a,b\}\in \pi }\frac{  \rho_{\eps}(\la a\, b\ra )}{ \eps N} \Delta(\pi).\end{aligned}  $$
Permuting its two sums, the first line of the right-hand-side can be rewritten as 
$$\sum_{1\le c<d\le n}\sum_{\substack{\pi \in\Pc\Mc(n)\\ \{c\},\{d\}\in \pi}} \esp\left[ \left(e^{\frac{\eps t}{2N}}U_t, \int_{0}^te^{\frac{ s\eps}{N}}U_s^{\ts 2}ds \right)^{\ts \pi}\right]\prod_{ \{i,j\}\in \check{\pi}^{c,d} }\frac{  \rho_{\eps}(\la i\, j\ra )}{ \eps N}.$$
Permuting again these sums and using (\ref{Explicit Projective family Casimir}), we  find that $\frac{d}{dt}R_t= R_t\Delta_G(n).$
\end{proof}
 
  The expectations on the right-hand-side of the above formulas are given by the Proposition \ref{BrownianSW}, yielding the Lemma \ref{exp unitary tensors} below. Let us define for any $1\le k\le n, $
 $$\Delta_\eps^n(k)=\frac{k}{\eps N}+ \Delta_\U(k)\ts \Id^{n-k}\in \End(V^{\ts n})$$
and for any $1\le k\le \min(n,m),$ $$\Delta^{n,m}(\min(n,m)-k)=\Delta_\U(n-k)\ts \Id^{\ts n+m-k}+\Id^{n+m-k}\ts \Delta_\U(m-k)\in\End(V^{n}\ts \overline{V}^{\ts m}).$$ We set $\Delta_\eps^n(0)=0$ and $\Delta^{n,m}(0)=0.$  We shall consider the set $\Sym[X_0,\ldots ,X_a]$ of symmetric  polynomials in $a+1$ variables $X_0,\ldots, X_a.$ For any $d\ge 1,$ we denote by  $H_d(X_0,\ldots,X_a)=\sum_{0\le i_1\le \ldots i_d\le a}X_{i_1}\ldots X_{i_d}=\sum_{\lambda_0,\ldots,\lambda_a\ge0: \lambda_0+\ldots+\lambda_a=d}X_0^{\lambda_0}\ldots X_{a}^{\lambda_a}\in \Sym[X_0,\ldots, X_a]$ the  complete symmetric polynomial of degree $d$ with the convention $H_0=1$, and we define an analytic function valued in symmetric polynomials setting for all $t\in \R,$
$$\Ec_t^{a,H}=\sum_{d\ge 0} \frac{t^{d+a}}{(d+a)!}H_{d}(X_0,\ldots, X_a).$$
When the context is clear, we shall drop the first upper index and denote this function by  $\Ec_t^H.$  Let us stress that these functions are not compatible: for all $a\ge 0,$  $\Ec_t^{a+1,H}(X_0,\ldots, X_a,0)\not=\Ec_t^{a,H}(X_0,\ldots, X_a).$ In return, they satisfy the following
\begin{lem}For all $t\in \R_+,$\label{Serie Gen Compl Simplex}
$$\Ec^{a,H}_t=\sum_{b=0}^a \frac{e^{tX_b}}{\prod_{c\not= b} (X_b-X_c)}=\int_{\Delta_{a}(t)}e^{t_0 X_0+\ldots + t_a X_a} d \mbf t , $$
where $\Delta_{a}(t)=\{ (t_0,\ldots,t_a)\in \R_+^{a+1}: t_0+\ldots +t_a=t\}$ and $d \mbf t$ denotes the Lebesgue measure on this simplex. 
\end{lem}
\begin{proof} Indeed, one easily checks that these three functions are the unique solution to the following problem: $(S^{a})_{a\ge 0}$ are   analytic functions, such that $S^0_t(X_0)= e^{tX_0}$ and for any $a\ge 1,$ $S^a\in \Sym[X_0,\ldots,X_a]$ satisfy for any  $t\ge0,$
\begin{equation}
\frac{d}{dt}S^a_t(X_0,\ldots, X_a)=X_a S^a_t(X_0,\ldots,X_a)+S^{a-1}_t(X_0,\ldots, X_{a-1}).\label{EDOSerieGen}
\end{equation}
\end{proof}
\begin{rema} Note that expanding the last column of the following, one can also rewrite \begin{equation*} \label{DetSerieExpCompl}
\Ec_t^H=\prod_{0\le i<j\le a}(X_j-X_i)^{-1}\det\left(\begin{array}{ccccc} e^{t X_0} & X_0^{a-1} & X_0^{a-2} &\cdots&  1 \\
 e^{t X_1} & X_1^{a-1} & X_1^{a-2} &\cdots&  1 \\
 \vdots& \vdots &  &\vdots & \vdots \\
e^{t X_a}& X_a^{a-1} & X_a^{a-2} &\cdots&   1
\end{array}\right).
\end{equation*}\end{rema}
For any $1\le a\le \frac{n}{2}$ and $1\le a'\le \min(n,m),$   let   $\mu_a\in \mathcal{P}\Mc(n)$ and $\nu_{a'}\in\Pc\Mc(n,m)$ be  the partial matchings with  respective $2$-blocks $\{n,n-1\},\ldots , \{n-2a+2,n-2a+1\} $ and $\{n,n+1\},\ldots, \{n-a'+1,n+a'\},$ whereas $\mu_0\in \Pc\Mc(n)$ and $\nu_0\in \Pc\Mc(n,m)$ are the partial matchings  with only blocks of size $1.$
When $\pi$ is a partial matching in respectively $\Pc\Mc(n)$  and $\Pc\Mc(n,n)$, we set $\Hc(\pi)=\{\sigma\in \Sy_{n}: \sigma (\mu_{|\pi|})=\pi\}$ and $\Dc(\pi)= \{\sigma\in \Sy_{n}\times \Sy_n: \sigma (\nu_{|\pi|})=\pi\}.$  
\begin{lem}\label{exp unitary tensors} The families $(\Delta_\eps^n(k))_{1\le k\le n}$ and $(\Delta^{n,m}(k))_{1\le k\le \min(n,m)}$ are  commutative.  For any $t\ge0$ and partial matching $\pi\in \Pc\Mc(n),$  
$$\esp\left[\left(e^{\frac{\eps t}{2N}}U_t, \int_{0}^te^{\frac{ s\eps}{N}}U_s^{\ts 2}ds \right)^{\ts \pi}\right]= \frac{2^{-|\pi|}}{(n-2|\pi|)!}\sum_{h\in \Hc(\pi)}h\Ec_t^H(\Delta^n_\eps(n),\Delta^n_\eps(n-2),\ldots, \Delta_\eps^{n}(n-2|\pi|))h^{-1}$$
and for any $\pi\in\Pc\Mc(n,m)$,
$$\begin{aligned}\esp &\left[\left((U_t, V_t),\int_{0}^tU_s\ts V_sds \right)^{\ts \pi}\right]\\&\hspace{2 cm}=\frac{1}{(n-|\pi|)!(m-|\pi|)!} \sum_{h\in \Dc(\pi)}h\Ec_t^H(\Delta^{n,m}(n),\Delta^{n,m}(n-1),\ldots, \Delta^{n,m}(n-|\pi|))h^{-1}.\end{aligned}$$
 \end{lem}
 \begin{proof} The first  property can been seen directly or  using the fact ii) of Lemma \ref{PropMB} stating that the law of the  Brownian motion is invariant by conjugation. Let us prove the formulas of the  Lemma by induction on $|\pi|$.  Again, we consider only the case $G\in\{\OO(N),\Sp(N)\}.$   When $|\pi|=0, $  according to Proposition  \ref{BrownianSW},  the left-hand-side is $\esp[(e^{\frac{\eps t}{2N}}U_t )^{\ts n}]= \exp\left[t\left(\Delta_U(n)+\frac{n\eps}{2N}\right)\right]=\exp\left(t\Delta^n_\eps(n)\right)=\Ec_t^H(\Delta^n_\eps(n)).$  For any $\mu\in \Pc\Mc(n), t\ge 0$, let us write $R_t(\pi)$ for the right-hand-side of the formula.  Let us use once more  the notation $\Delta(\pi)\in\End(V^{\ts n})$ that has been defined in (\ref{NotLaplacienAppar}). Note that for any $h\in \Hc(\pi),$ $\Delta(\pi)= h \Delta^n_\eps(n-2|\pi|).$   Using (\ref{EDOSerieGen}),  we get  
 $$\begin{aligned}&\frac{d}{dt}R_t(\pi)= R_t (\pi)\Delta(\pi) +\sum_{\{a,b\}\in\pi}R_t(\hat{\pi}_{a,b}),\end{aligned}$$
where, as above, $\hat{\pi}_{a,b}$ denotes the matching with $\{a,b\}$ replaced by $\{a\},\{b\},$ when $\{a,b\}\in\pi.$

Let us recall that according to the equation (\ref{Diff Matched Tensors}),  the family  $\{\esp\left[\left(e^{\frac{\eps t}{2N}}U_t, \int_{0}^te^{\frac{ s\eps}{N}}U_s^{\ts 2}ds \right)^{\ts \pi}\right]:\pi\in \Pc\Mc(n) \}$ satisfies the same system of ordinary differential equations. Using the induction assumption yields then the announced equality.
 \end{proof}
 
 We want now to  deduce from the Theorem  \ref{Unitary representation} an explicit expression for the expectations of  monomials in entries of $G_t$ in terms of some functions on $\Bc_n(\eps N)$ or $\Bc_{n,m}(N).$  

\vspace{0,2 cm}

\noindent\emph{Orbits  under $\Sy_n$ left and right actions on $\Bc_n$:} Let us consider the action of  $\Sy_n\times \Sy_n$ on $\Mc(n)$ given for any $\alpha,\beta\in \Sy_n$ and $\pi\in \Mc(n)$  by $(\alpha,\beta).\pi=\alpha\pi\beta^{-1}.$ For any element $\pi\in \Bc_n,$ let us write $h(\pi)$ for half  the \emph{number of horizontal curves} of the diagram of $\pi.$ Two matchings $\mu,\nu\in\Mc(n)$ are in the same orbit if and only if $h(\mu)=h(\nu).$  Let us highlight representatives of these orbits setting for any integer $1\le a\le \frac{n}{2}$,
$$\pi_a=\prod_{k=0}^{a-1}\la n-2k-1\; n-2k\ra.  $$
Let us recall that $\pi_0$ also denotes the identity of $\Bc_n(z).$ Similarly $\left(\Sy_n\times \Sy_m\right)^{\times 2}$ acts on $\Mc(n,m),$  the function $h$ characterizes its orbits and we set for any $1\le a\le \min(n,m),$ 
$$\tilde{\pi}_a= \prod_{k=0}^{a-1} \la n-k \;  n+k+1\ra. $$

\vspace{0,2 cm}
\noindent\emph{Deformed Weingarten function:}   Let us fix a notation for the sum of permutation like  elements appearing in Proposition \ref{BrownianSW}, setting for any integer $l\in \{1,\ldots,n\},$ 
$$-Z_l= \frac{l(1-\eps/N)}{2}+\frac{1}{ N}\sum_{1\le i<j\le l}(i\, j)\in \C[\Sy_n]$$
and $Z_0=0.$
We shall also set for $l\le \min(n,m),$
$$-Y_{\min(n,m)-l}=\min(n,m)-l+\frac{1}{N}\sum_{1\le i<j\le n-l }(i\,j)+\frac{1}{N}\sum_{ n+l\le i<j\le n+m }(i\,j)    \in \C[\Sy_n\times \Sy_m].$$
Notice that these two families of operators are commutative.   We set for all $t\ge 0,$ and all integer $0\le a \le \frac{n}{2},$
$$W^a_t= \frac{1}{(2\eps N)^a(n-2a)!}\Ec_t^H(Z_n,Z_{n-2},\ldots, Z_{n-2a})$$
and for $a\le \min(n,m),$
$$\mathcal{W}^a_t= \frac{1}{ N^a(n-a)!(m-a)!}\Ec_t^H(Y_n,Y_{n-1},\ldots, Y_{n-a}).$$
\noindent \emph{Kernel and sign of indices functions:} For any $I=(i_k)_{1\le k\le 2n}\in \{1,\ldots, N\}^{2n},$    let  $\ker_{1}(I)$  (or simply $\ker(I)$) and $\widetilde{\ker}_{-1}(I)$ be  the set of matchings $\pi\in \Mc(n),$ such that  for all  blocks $\{a,b\}$ of $\pi,$   we have  respectively $i_a=i_b$  and  $i_a=i_b+\frac{N}{2} \mod(N).$   For any permutation $\sigma\in\Sy_{2n},$ let us consider $${\iv}(\sigma,I)={\#\{k\in\{1,\ldots,p\}: i_{\sigma(2n-k+1)}<i_{\sigma(k)}\}}.$$  If $\pi\in \ker_{-1}(I),$  the quantity   $\eps(\sigma) (-1)^{\iv(\sigma^{-1},I)}$ does not depend on the choice of $\sigma\in \Sy_{2n}$ such that $\sigma(\pi_0)=\pi.$   For any indices $(i_k)_{1\le k\le 2n}\in \{1,\ldots, N\}^{2n},$ let us define $I^\pi=(i'_k)_{1\le k\le 2n}\in \{1,\ldots, N\}^N$ setting $i'_k= i_k+\frac{N}{2}\mod( N),$ for all $n< k\le 2n$ that is matched by $\pi$ with an element of $\{1,\ldots,n\}$,  and $i'_k=i_k$ otherwise. We then consider $\ker_{-1}(\pi)=\{I^{\pi}: I\in \widetilde{\ker}_{-1}(\pi)\}$ and set for any $I\in \ker_{-1}(\pi),$ 
$$\eps_I(\pi)= \eps(\sigma) (-1)^{\iv(\sigma,I^{\pi})},$$
for any $\sigma\in \Sy_{2n}$ such that $\sigma(\pi_0)=\pi.$ Let us also denote by $\la\cdot,\cdot\ra$  the inner product on $V^{\ts n}$ such that for all $v_1,\ldots,v_n\in V,$ 
\begin{equation}
\|v_1\ts\ldots \ts v_n\|^2= \|v_1\|^2\ldots \|v_n\|^2.\label{DefIP}
\end{equation}
Formulas (\ref{RepO}) and (\ref{RepSp}) show that the above definitions are tailored so that for any $\pi\in\Mc(n)$ and $I\in \{1,\ldots,N\}^{2n},$ 
$$\la e_{i_{1}}\ts \ldots \ts e_{i_{n}},\rho_{\eps}(\pi)e_{i_{2n}}\ts \ldots \ts e_{i_{n+1}}\ra =\eps_I(\pi)^{(1-\eps)/2},$$
if $I\in \ker_\eps(\pi)$ and $0$ otherwise.
     We can now give an equivalent statement  to Theorem \ref{Unitary representation}.  

\vspace{0,2 cm}

\begin{theo}\label{Deformed Weingarten}For any $I=(i_k)_{1\le k\le 2n}\in \{1,\ldots, N\}^{2n},$
$$\esp[(G_t)_{i_1,i_{2n}}\ldots (G_t)_{i_{n-1},i_{n+2}}(G_t)_{i_n,i_{n+1}}]=\sum_{\mu\in \ker_\eps(I)}\eps_I(\pi)^{(1-\eps)/2}\sum_{\alpha,\beta\in\Sy_n: \,\alpha\pi_{h(\mu)}\beta=\mu}W^{h(\mu)}_{t}(\beta\alpha), $$
 if $G=\OO(N)$ or $\Sp(N)$. When $G=\U(N),$
 $$\begin{aligned}\esp[(U_t)_{i_1,i_{2(n+m)}}\ldots(U_t)_{i_n,i_{n+2m+1}}&(\overline{U}_t)_{i_{n+1},i_{n+2m}}\ldots (\overline{U}_t)_{i_{n+m},i_{n+m+1}}]\\&=\sum_{\mu\in \ker(I)}\sum_{\alpha,\beta\in \Sy_n\times \Sy_m:\,\alpha \tilde{\pi}_{h(\mu)}\beta=\mu}\Wc^{h(\mu)}_{t}(\beta\alpha). \end{aligned}$$
\end{theo}

\begin{proof}[Proof of Theorem \ref{Deformed Weingarten}]This is a direct application of  Theorem \ref{Unitary representation} together with Lemma \ref{exp unitary tensors}. 
\end{proof}

\begin{rema} For any $\zeta\in \C,$ $\alpha\in \{0,1\},$ let us consider  
$$-\Delta^\alpha_\zeta(n)=\frac{n}{2}(\zeta-1)+\sum_{1\le i<j\le n}((i\, j)  -\alpha\la i\,j\ra    )\in \Bc_n(\zeta)$$
and 
$$-\Delta^\alpha_\zeta(n,m)=\frac{n\zeta}{2}+\sum_{1\le i<j\le n+m }(\un_{\{i,j\le n \text{ or }i,j>n\}}(i\, j)  -\alpha\un_{\{i\le n<j\}}\la i\,j\ra    )\in \Bc_{n,m}(\zeta).$$
 For any $n'\ge n,$ $m'\ge m,$ let us consider the embedding of the algebras $\Bc_{n}(\zeta),$ $ \Bc_{n,m}(\zeta),$  respectively into  $\Bc_{n'}(\zeta)$ and $ \Bc_{n',m'}(\zeta),$ such that any diagram $\pi \in \Bc_n$ (resp.   $\Bc_{n,m}$) is sent to a diagram in resp.  $\Bc_{n'}$ and $\Bc_{n',m'},$ having vertical lines going through points labeled by $\{n+1,\ldots, n'\}$ (resp. $\{n+1,\ldots,n'\}\times \{n'+1,\ldots, n'+m'-m\}$).  Note that 
 $$\la n+1\, n+2\ra \Delta_{\zeta}^{1}(n+2)=\la n+1\,\, n+2\ra\Delta_{\zeta}^{1}(n) \in \Bc_{n+1}(\zeta)$$
 and 
 $$\la n\, n+1\ra \Delta_{\zeta}^{1}(n+1,m+1)=\la n\, \,n+1\ra\Delta_{\zeta}^{1}(n,m) \in \Bc_{n+1,m+1}(\zeta).$$
 These equalities are generalization of (\ref{Projective family Casimir}). Following the above proof, it can also be shown that for any $z,\zeta\in \C,$    
 $$\begin{aligned}\exp( z \Delta^1_\zeta(n))= \sum_{\pi\in\Pc\Mc(n)} W_z^\zeta(\pi) \prod_{ \{a,b\}\in \pi } \la a\, b\ra \in \Bc_{n}(\zeta)
\end{aligned}$$
 and 
 $$\begin{aligned}\exp( z \Delta^{1}_\zeta(n,m))= \sum_{\pi\in\Pc\Mc(n,m)} \Wc_z^\zeta(\pi) \prod_{ \{a,b\}\in \pi }\la a\, b\ra \in \Bc_{n,m}(\zeta),
\end{aligned}$$
where $$W_z^\zeta(\pi)=\frac{|\pi|!}{|\Hc(\pi)|}\sum_{h\in \Hc(\pi)}h\Ec_z^H(\Delta^0_\zeta(n),\Delta^0_\zeta(n-2),\ldots, \Delta^0_\zeta(n-2|\pi|))h^{-1}\in \C[\Sy_n]$$
and  $\Wc_z^\zeta(\pi)$ is 
$$ \frac{|\pi|!}{|\Dc(\pi)|} \sum_{h\in \Dc(\pi)}h\Ec_z^H(\Delta^{0}_\zeta(n,m),\Delta^0_\zeta(n-1,m-1),\ldots, \Delta^{0}_\zeta(n-|\pi|,m-|\pi|)) h^{-1}\in \C[\Sy_{n+m}].$$
\end{rema}

\section{Integration against the Haar measure and FFT}

Let us refer here the Reader to \cite{ThCollins,CS,CM,ZJ} where he can find  a classical approach to integration over classical groups. These approaches use orthogonality of irreducible characters together with the decomposition of the tensor space into  irreducible representations or the first fundamental theorem of invariants. The latter is used in the following way (see \cite{CM,ZJ}): by left-invariance of the Haar measure, the expectation of the representation of a Haar distributed random variable is equal  to the projection on the invariant space. Besides, the FFT theorem gives a generating family of this space. To compute the projection, and thereby the mean of the representation, one way is therefore to find a pseudo inverse for the Gram matrix of this family.

\begin{figure}[!h] 
  \centering
 \includegraphics[height=0.8 in, width= 1.5 in]{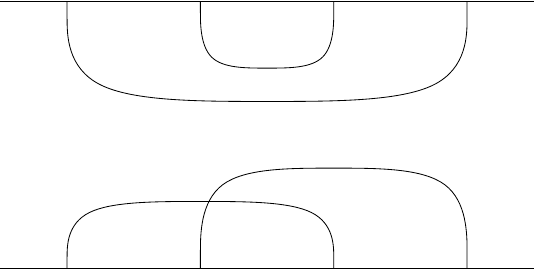}
  \caption{Diagram $\pi_{\mu,\mu'}= \{\{1,4\},\{2,3\},\{5,7\},\{6,8\}\}$, where $\mu=\{\{1,4\},\{2,3\}\}$ and $\mu'=\{\{1,3\},\{2,4\}\}$ .}\label{JuxtapoMatching}
\end{figure}

\vspace{0,5 cm}
Let us give here another approach considering the behavior of the expectation with respect to Brownian motions as $t\to\infty.$ Let $(U_t)_{t\ge 0},$  $(V_t)_{t\ge 0}$ be two $\U(N)$-Brownian motions and $H$ be a Haar distributed random variable on $G$. Let  $\eps_N= \frac{1}{N!}\sum_{\sigma\in \Sy_N}\eps(\sigma)\sigma$ and recall that  $\la\cdot,\cdot\ra$ denotes the inner product on $V^{\ts n }$ such that (\ref{DefIP}) holds.
\begin{lem} \label{Spectrum} The endomorphisms $(\Delta_\eps^n(k))_{1\le k\le n},$  $(\Delta^{n,m}(k))_{1\le k\le \min(n,m)}$  are self-adjoint operators with   nonpositive spectrum and are all  invertible but $\Delta_1^n(N),$  that satisfies $\ker(\Delta_{1}^n(N))= \Im(\rho(\eps_N)).$ The following holds as $t\to\infty.$ Let $\pi\in\Pc\Mc(n)$  be a partial matching with $|\pi|<\frac{n}{2}.$  If  $n\not=2|\pi|+N$ or $\eps=-1, $ then  $$\esp\left[\left(e^{\frac{\eps t}{2N}}U_t, \int_{0}^te^{\frac{ s\eps}{N}}U_s^{\ts 2}ds \right)^{\ts \pi}\right]\to0.$$
If $n=2|\pi|+N,$ then for any $\sigma\in \Sy_n$ sending the $2$-blocks of $\pi $ to $\{N+1,N+2\},\ldots, \{n-1,n\},$ $$\esp\left[\left(e^{\frac{ t}{2N}}U_t, \int_{0}^te^{\frac{ s}{N}}U_s^{\ts 2}ds \right)^{\ts \pi}\right]\rho\left(1-\sigma^{-1}\eps_N\sigma\right)\to0 .$$
If $\pi\in\Pc\Mc(n,m),$ with $|\pi|<\min(n,m),$ then 
 $$\esp\left[ \left((U_t, V_t),\int_{0}^tU_s\ts V_sds \right)^{\ts \pi}  \right] \to 0.$$
\end{lem}
\begin{proof} Any endomorphism $\rho_\eps( \tau)\in \End(V^{\ts n}),$  with   $\tau\in \Sy_n$  a transposition,  is self-adjoint. In particular, all the operators considered are self-adjoint. For any $\lambda_0,\ldots,\lambda_a\in \left(\R_-^{*}\right)^{a+1},$ $\Ec_t^H(\lambda_0,\ldots,\lambda_a)=\int_{\Delta_a(t)}\prod_{k=0}^a e^{\lambda_k t_k}d\mbf t\to0$ and using Lemma \ref{exp unitary tensors}  and \ref{Serie Gen Compl Simplex}, it is enough to prove that all the operators considered have a negative spectrum.   Using the expression of Proposition \ref{BrownianSW} and bounding the operator norm of $\rho_\eps( \tau)$ by $1$ for any transposition $\tau$ yield that the operators  $(\Delta_\eps^n(k))_{ k<N},$  $(\Delta^{n,m}(k))_{1\le k\le N}$ and $\Delta^n_{-1}(N)$ have  negative spectrum,  whereas $\Delta^n_{1}(N)$ have non-positive spectrum.  An elementary inspection leads then to  $\ker(\Delta_{1}^n(N))= \Im(\rho(\eps_N)).$  Let $(G_t)_{t\ge 0}$ be a Brownian motion on $\SU(N)$ and let $\Delta_\SU(k)\in\End(V^{\ts k})$ be such that $\esp[S_t^{\ts k}] = \exp(t\Delta_{\SU}(k)).$   According to Lemma \ref{Facto U SU Cercle},  $\Delta_\U(k)= \Delta_\SU(k)- \frac{k^2}{2N^2}.$  The operator $\Delta_{\SU}(k) $ is self-adjoint and as $\esp[S_t^{\ts k}] $ is a contraction, it has nonpositive spectrum.  Therefore,  when $k >N,$   
$$\Delta^n_\eps(k)=- \frac{k^2-\eps k N}{2N^2} +\Delta_{\SU}(k)\ts \Id^{n-k} $$
and setting $k'= \min (n,m)-k,$
$$\begin{aligned}\Delta^{n,m}(k)=- \frac{(n-k')^2+(m-k')^2}{2N^2}+ \Delta_\SU(n-k')\ts \Id^{\ts n+m-k'}+\Id^{n+m-k'}\ts \Delta_\SU(m-k')\end{aligned}$$
have negative spectrum.
\end{proof}
\begin{proof}[Proof of Theorem \ref{FFT}]  Let us recall  from section  \ref{SectionFFT} that it is equivalent to prove that the vector spaces $\left(V^{\ts 2n}\right)^{\OO(N)}, $     
$\left(V^{\ts 2n}\right)^{\Sp(N)}$ and $\left(V^{\ts n}\ts \overline{V}^{\ts n}\right)^{\U(N)}$ are spanned  respectively by   $(w_\pi)_{\pi \in\Mc(n)},$  respectively  for $\eps\in\{1,-1\},$ and $(w_\pi)_{\pi \in \Mc_+(n)},$  as defined in (\ref{Spanning Family of Invariant}) and (\ref{signSymp}).   We denote in any each case by $\mathcal{I}_G$ the  linear span of the latter families.  For any unitary representation  $(\rho,W)$ of $G,$  denoting by $H$ a Haar distributed random variable on $G,$ $\esp[ \rho( H)]$ is equal to the Hermitian projection on $W^G.$  Let us consider a $G$-Brownian motion $(G_t)_{t\ge0},$ respectively  issued from $\Id$, when $G$  is $\U(N)$ or $\Sp(N)$  and from a random variable $S$ such that $\E[\det(S)]=0$, when $G=\OO(N).$  When  $W$ is     
$V^{\ts 2n}$ or  $V^{\ts n}\ts \overline{V}^{\ts n}$ and    $G$  is $\OO(N),\Sp(N)$ or $\U(N), $  according  to Theorem \ref{Unitary representation} and Lemma \ref{Spectrum},  for any $t\ge 0,$ $\esp[\rho(G_t)]=e^{t\Delta},$ where $\Delta$ is Hermitian and nonpositive and as $t\to \infty,$   $\esp[\rho(G_t)]\to P,$ where $\Im P\subset \Ic_G.$ The endomorphism $P$ is  an Hermitian projection. Moreover, by left-invariance of the Haar measure, $P\esp[\rho(H)]= \esp[\rho(H)].$ Hence, $W^G\subset \Ic_G$ (and  $P=\esp[\rho(H)]$).
\end{proof}
Note that the latter argument gives an  explicit way to prove that $\esp[G_t^{\ts n}]\to \esp[H^{\ts n}],$ as $t\to \infty$,  where $H$ is distributed according to the Haar measure on $G$, and $(G_t)_{t\ge 0}$ is a $G$-Brownian motion with a proper initial condition. It therefore leads to a new way to compute moments for the Haar measure.  For any partial matching in respectively $\Pc\Mc(n)$  and $\Pc\Mc(n,n)$, let us set $\Hc(\pi)=\{\sigma\in \Sy_{n}: \sigma (\mu_n)=\pi\}$ and $\Dc(\pi)= \{\sigma\in \Sy_{n}\times \Sy_n: \sigma (\nu_n)=\pi\}.$  
\begin{lem} \label{formula Haar from Brownian}  Let us denote  by $H$ a Haar distributed random variable on $G $ and by $(G_t)_{t\ge 0}$ a Brownian motion on $G,$ respectively  issued from $\Id$, when $G$  is $\Sp(N)$  and from a random variable $S$ such that $\E[\det(S)]=0$, when $G=\OO(N).$ Then,    
$$\begin{aligned}&\esp[G_t^{\ts n}]\to \esp[H^{\ts n}]\\&=\sum_{\pi\in\Pc\Mc(n): 2|\pi|=n}\left(\frac{|\pi|!}{|\Hc(\pi)|}\sum_{h\in \Hc(\pi)}h(\Delta_\eps^n(n)\Delta_\eps^n(n-2)\ldots \Delta_\eps^n(2))^{-1}h^{-1}\right)\prod_{ \{a,b\}\in \pi }\frac{  \rho_{\eps}(\la a\, b\ra )}{- \eps N}.
\end{aligned} $$
When $G=\U(N),$
$$\begin{aligned}&\esp\left[U_t^{\ts n}\ts \overline{U}_t^{\ts n}\right]\to \esp\left[H^{\ts n}\ts \overline{H}^{\ts n}\right]\\
&=\sum_{\pi\in \Pc\Mc(n,n): |\pi|=n} \left(\frac{|\pi|!}{|\Dc(\pi)|}\sum_{h\in \Dc(\pi)}h(\Delta^{n,n}(n)\Delta^{n,n}(n-1)\ldots \Delta^{n,n}(1))^{-1}h^{-1}\right) \prod_{ \{a,b\}\in \pi }  \frac{ \rho_{+}(\la a\, b\ra )}{-N}.\end{aligned}$$
\end{lem}
\begin{proof} As argued in the above proof of Theorem \ref{FFT}, as $t\to \infty, $  $\esp[G_t^{\ts n}]\to \esp[H^{\ts n}]$. Using Theorem \ref{Unitary representation}, Lemma \ref{exp unitary tensors} and the first expression  of Lemma \ref{Serie Gen Compl Simplex} together with the first part of Lemma \ref{Spectrum} yields the statement.
\end{proof}
The latter expressions can be simplified as follows. When $W$ is a subgroup of $\Sy_n,$ we shall consider the idempotent $\Pi_{W}= \frac{1}{|W|}\sum_{w\in W}w$ and set $\Pi^{\eps N}_W=\rho_\eps(\Pi_W).$  Let us introduce for any integer  $1\le a\le n,$
$$X_a=\sum_{b<a} (b\, a)\in \C[\Sy_{n}].$$
These elements form a commutative family  and are called \emph{Jucys-Murphy elements.} These elements enjoy a lot of properties (see \cite{Jucys,OkounkovVershik,OkounkovApproachItal}), we shall recall two of them (see Proposition 1 and 3 of \cite{ZJ} for a proof). For any integer $k\in\N^*,$  let us consider the \emph{lattice of  partitions} of $\{1,\ldots,k\}$ given by inclusion.   For any pair of partitions $\pi$ and $\nu$,   $\pi\vee \nu$ is the finest partition coarser than $\pi$ and $\nu$. When $\pi$ is respectively a partition  or a permutation, we denote by $\#\pi$  its \emph{number of blocks or its number of cycle}.  This function is related to Jucys-Murphy elements as follows:  for any $\zeta\in\C,$
\begin{equation}
\prod_{k=1}^n (\zeta+X_{2k-1}) \pi_0= \sum_{\mu\in \Mc(n)}  \zeta^{\# \mu\vee \mu_{n} } \mu\label{FactoDemiOmega}
\end{equation}
and 
\begin{equation}
\prod_{k=1}^n (\zeta+X_{k}) = \sum_{\sigma\in \Sy_n} \zeta^{\#\sigma}\sigma.\label{FactoOmega}
\end{equation}
Note that these  formulas imply that $\prod_{k=1}^n (\zeta+X_{2k-1}) $ commutes with $\Pi_{\Hc_{n}}$ and that $\prod_{k=1}^n(\zeta+X_k)$ is in the center of $\C[\Sy_n].$ 
\begin{lem}\label{Reduction} For any integer $n\ge 1,$
$$(-\eps N)^nn!\rho_\eps\left(\prod_{1\le a\le n } (\eps N+X_{2a-1})\right)\Pi^{\eps N}_{\Hc_{n}}(\Delta_\eps^{2n}(2n)\Delta_\eps^{2n}(2n-2)\ldots \Delta_\eps^{2n}(2))^{-1}\Pi^{\eps N}_{\Hc_{n}}=\Pi^{\eps N}_{\Hc_{n}}$$
and 
$$(-N)^nn!\rho\left(\prod_{1\le a\le n} ( N+X_{a})\right)\Pi^N_{\Dc_n}(\Delta^{n,n}(n)\Delta^{n,n}(n-1)\ldots \Delta^{n,n}(1))^{-1}\Pi_{\Dc_n}^N=\Pi_{\Dc_n}^N.$$
\end{lem}
\begin{proof}   First notice that \begin{equation*}
\sum_{1\le a<b\le 2n}(a\,b)\Pi_{\Hc_n} =\Pi_{\Hc_n}n(1+X_{2n-1})\Pi_{\Hc_n}\tag{*}\label{Quotiented Casimir}
\end{equation*}

and 
$$\sum_{\substack{1\le a<b\le n\\ \text{or } n<a<b\le 2n}} ( a\hspace{0,2cm} b) \Pi_{\Dc_n}= \Pi_{\Dc_n}n X_n\Pi_{\Dc_n}.$$
Indeed, on the right-hand-side  $\Pi_{\Hc_n}X_{2n-1}\Pi_{\Hc_n}= 2(n-1)\Pi_{\Hc_n}\left(2n-2\hspace{0,2cm} 2n-1\right)\Pi_{\Hc_n}  $ and on the left-hand-side,  each transposition  $( 2k-1\hspace{0,2cm} 2k)$  with $1\le k\le n$ belongs to $\Hc_n,$ whereas   the $2(n-1)$ others  are conjugated to  $(2n-2 \hspace{0,2cm} 2n-1)$ by an element of $\Hc_n$. As $\sum_{1\le a<b\le 2n} (a\, b)$ is  central in $\Sy_{2n}$ the first formula holds true.  A similar argument with respectively  $\Dc_n$  and $\Sy_n\times \Sy_n$ in place of $\Hc_n$  and   $\Sy_{2n}$ yields the second one.      Let us set $\Pi_n=\Pi_{\Hc_n}^{\eps N}$, and for any $1\le k\le n$, $F_k=-\eps N\Delta_\eps^{2n}(2k)$ and $L_k=k \rho_\eps\left(\eps N+X_{2k-1}\right)$  (resp.  $\rho(\Pi_{\Dc_n}),$ $N\Delta^{n,n}(k)$ and $k\rho\left(N+X_k\right)$).  We consider $\Rc_n=\prod_{k=1}^n F_k^{-1},$ and $\Omega_n=\prod_{k=1}^nL_k$.   The two equalities above imply that 
\begin{equation*}
F_n=\Pi_n L_n\Pi_n, \tag{**}\label{Relation Casimir JucysM}
\end{equation*} for any integer $n\ge 1.$  The claim is equivalent to the following equality
\begin{equation}
\Pi_{n} \Omega_n \Rc_n\Pi_n=\Pi_n \label{reduction inversion},
\end{equation}
for all $n\ge 0,$ that we shall prove by induction. Consider as initial step $n=0,$ where  factors  are empty and equal to $1$, and the statement is trivial.
Assuming the equality for $n\ge 0$ and then using (\ref{Relation Casimir JucysM})  yields
$$\begin{aligned}\Pi_{n+1}\Omega_{n+1}\Rc_{n+1}\Pi_{n+1}&= \Pi_{n+1}L_{n+1}\Omega_{n} \Rc_{n} F_{n+1}^{-1} \Pi_{n+1}\\
&= \Pi_{n+1}L_{n+1}\Pi_n\Omega_{n} \Rc_{n}\Pi_n  F_{n+1}^{-1} \Pi_{n+1}\\
&= \Pi_{n+1}L_{n+1}\Pi_n  F_{n+1}^{-1} \Pi_{n+1}= \Pi_{n+1}L_{n+1}\Pi_{n+1}F_{n+1}^{-1}=\Pi_{n+1},\end{aligned}$$
where we have furthermore used  that $F_{n+1}$ commutes with the action of $\Sy_{2n+2}$ (resp. $\Sy_{n+1}\times\Sy_{n+1}$), $L_{n+1}$ with $\Pi_{n}.$ \end{proof}
Let us consider respectively  an element $W_n^\eps$ of $\C[\Sy_{2n}]$ invariant by left or right translation by $\Hc_n$ and   a central element $\Wg_n$ of  $\C[\Sy_{2n}],$ such that 
$$   \rho_\eps(W^{\eps}_n)=(-\eps N)^nn! \Pi^{\eps N}_{\Hc_n}(\Delta_\eps^{2n}(2n)\Delta_\eps^{2n}(2n-2)\ldots \Delta_\eps^n(2))^{-1}\Pi^{\eps N}_{\Hc_n} $$
and 
$$\sum_{\alpha,\beta} \Wg_n(\alpha\beta^{-1})\rho(\alpha)\ts\rho(\beta) =(- N)^nn!\Pi^N_{\Dc_n}(\Delta^{n,n}(n)\Delta^{n,n}(n-1)\ldots \Delta^{n,n}(1))^{-1}\Pi_{\Dc_n}^N.$$
The Lemma  \ref{Reduction} yields that $\rho_\eps(W_n)$ and $\sum_{\alpha,\beta} \Wg_n(\alpha\beta^{-1})\rho(\alpha)\ts\rho(\beta)$ are respectively the \emph{pseudo-inverses} of  the two self-adjoint operators    $\sum_{\sigma\in\Sy_{2n}}  N^{\# \sigma\mu_n \vee \mu_n }\rho_{\eps}(\sigma )$  and $\sum_{\alpha,\beta\in\Sy_{n}}N^{\#\alpha \beta^{-1}} \rho(\alpha\times \beta).$ 
We can reformulate Lemma \ref{formula Haar from Brownian} as follows. For any pair of matchings $(\mu,\mu')\in\Mc(n),$ let us consider  the matching $\pi_{\mu,\mu'}\in \Mc(2n)$  whose diagram is given by $\mu$ in the positive half plane and  by  $\mu'$ in the negative half plane (see figure \ref{JuxtapoMatching}). Then,  Lemma \ref{formula Haar from Brownian} can be reformulated into

\begin{theo}[\cite{CS,CM,ZJ}] If $G$ is $\OO(N)$ or $\Sp(N)$, then
$$ \esp[H^{\ts 2n}]= \sum_{\alpha,\beta\in \Sy_{2n}} W^{\eps}_n(\alpha\beta^{-1}) \rho_\eps\left(\pi_{\alpha\mu_n,\beta \mu_n}\right).$$
If $G=\U(N),$ 
$$\esp[H^{\ts n}\ts \overline{H}^{\ts n}]= \sum_{\alpha,\beta\in \Sy_{n}} \Wg_n(\alpha\beta^{-1}) \rho_+(\pi_{\alpha \nu_n, \beta \nu_n}). $$
\end{theo}
The two  pseudo-inverses appearing in these formulae can be re-expressed using representation theory. For example, it can be shown that for any $k\in\{1,\ldots, 2n\},$ $\rho_{\eps}(\eps N+ X_{k})$  is invertible in $\End(V^{\ts 2n}),$ so that  
$$ \rho_\eps(W^{\eps}_n)= \rho\left( \prod_{k=1}^n(\eps N+ X_{2k-1})\right)^{-1} \rho_\eps(\Pi_{\Hc_n})$$
and 
$$\sum_{\alpha,\beta\in \Sy_n} \Wg_n(\alpha\beta^{-1})\rho(\alpha)\ts\rho(\beta)= \rho\left(\prod_{k=1}^n(N+X_k)\right)^{-1} \rho\left(\Pi_{\Dc_n}\right).$$
Using these two equalities together with (\ref{FactoDemiOmega}) and $(\ref{FactoOmega}),$ leads to the exact statement  of the formulae of \cite{ZJ},  that is equivalent to  the one of \cite{CS}.	

\begin{rema}The above way to get formulae for the Haar measure can also be applied to $\SU(N).$  This case has not been  explored as much as the three series $\OO(N),\Sp(N)$ and $\U(N)$.  Let us  consider an example. If $S$ is Haar distributed on $\SU(N)$, then multiplying on the left by rotation matrices leads to 
$$\esp[ \prod_{1\le i,j\le N} S_{i,j}]=0,$$
when $N$ is odd.  On the  other hand, the following is a conjecture  \cite{Landsberg}:  that  for all $N$ even,
\begin{equation}
\esp[  \prod_{1\le i,j\le N} S_{i,j}]\not= 0. \label{Alon Tarsi}
\end{equation}
This  conjecture has been shown in \cite{Landsberg} to be equivalent to the Alon-Tarsi conjecture (\cite{AT}) about even and odd Latin squares, as well as Hadamard-Howe and Foulkes conjectures (see \cite{Landsberg} for references).
In \cite{Glynn,Drisko}, it has been  proved to hold true for any even number such that $N\pm1$ is a prime number. Though, the general case is still open.  Using $\SU(N)$-Brownian motion, it can be reformulated as follows.  Let $(S_t)_{t\ge 0}$ be a $\SU(N)$-Brownian motion  issued from identity.   As $\esp[S_t^{\ts N^2}]\to \esp[S^{\ts N^2}],$ when $t\to\infty, $ according to Lemma \ref{Facto U SU Cercle} and  Proposition \ref{BrownianSW}, the Alon-Tarsi conjecture is equivalent to the following: for any even integer $N,$  
\begin{equation}
 \la e_1^{\ts N}\ts \ldots \ts e_{N}^{\ts N}, \exp\left(-t \sum_{1\le i<j\le N^2}(i\, j)\right) (e_1\ts \ldots \ts e_N)^{\ts N} \ra\to x_N\not =0,\label{Brownian Alon Tarsi}
\end{equation}
as $t\to\infty.$
\end{rema}

\appendix
\section{Symplectic invariants\label{Proof Sp}}
\begin{proof}[Proof of Lemma \ref{algebra morphisms} for $\Bc_n(-N)$]  As $\{(a\, b), \la a\, b\ra : 1\le a<b\le n\}$ generates the algebra $\Bc_n(-N)$, it is sufficient to prove that for any $\pi\in \Mc(n),$ $1\le a<b\le n,$ $\rho_{-}(\pi (a\,b))=\rho_-(\pi)\rho_-((a\, b))$ and  $\rho_{-}(\pi \la a\,b\ra )=\rho_-(\pi)\rho_-(\la a\, b\ra).$  On  the one hand, $\pi .(a\, b)= (2n+1- a\, 2n+1-b) \left(\pi\right), $ where the left-hand-side is a product in $\Bc_n(-N)$ and the right-hand-side an action of $\Sy_n$ on $\Mc(n).$ On the other hand, if $\pi$ matches $2n+1-a$ and $2n+1-b$ with  respectively   $a'$ and $b',$ then whether $a'=2n+1-b$ and $b'=2n+1-a,$ in which case, $\pi \la a \,b\ra =-N\pi $, or $\pi \la a\, b\ra=(2n+1-a\,\, b') \left(\pi\right) =(2n+1-b\,\, a')\left(\pi\right) $.  Using  (\ref{RepSp}) and (\ref{RepContracSp}), it is now elementary to check the multiplicativity in both cases.
\end{proof}

\begin{proof}[Proof of proposition \ref{BrownianSW} for $G=\Sp(N)$] Let us define  $\iota:M_{N/2}(\C)\times M_{N/2}(\C)\to M_{N}(\C): (A,B)\mapsto \left(\begin{array}{cc} A & -\overline{B} \\ B  & \overline{A}\end{array}\right) $ and recall that $$\usp=\{\iota(A,B) : A\in \uN, B\in M_N(\C), {}^tB=B \}.$$ Let us choose the basis  of $\usp$ formed by the unions of the following families: 

$$\frac{1}{\sqrt{2N}}\{\iota( E_{a,b}-E_{b,a},0),\iota(0, E_{a,b}-E_{b,a}),\iota(i( E_{a,b}+E_{b,a}),0),\iota(0,i( E_{a,b}+E_{b,a})):  1\le a<b\le \frac{N}{2}  \}$$  and 

$$\frac{1}{\sqrt{N}}\{\iota(i E_{a,a},0), \iota(0, E_{a,a}),\iota(0, iE_{a,a}) :1\le a\le \frac{N}{2}\}.$$ 
The Casimir element of $\usp$, that is $c_\usp=\sum_{1\le i\le \frac{N^2}{2}}x_i\ts x_i,$ where $(x_i)_{1\le i\le N^2 }$ is the above orthogonal basis, viewed as an element of the complex enveloping algebra $\Uc(\glN(\C))$, has the following expression:  
\begin{equation*}-Nc_\usp=-\sum_{1\le a,b,c,d\le N}J_{a,c}J_{b,d}E_{a,b}\ts E_{c,d}+\sum_{1\le a,b\le N}E_{a,b}\ts E_{b,a}\in \Uc(\glN(\C)). \end{equation*}
Considering  $\rho_{V^{\ts n}}$ as a  representation of the algebra $\Uc(\glN(\C))$,  

\begin{align*}
-N\rho_{V^{\ts n}}(\Delta_{\Sp(N)})&=2(N+1)n+ 2\sum_{1\le a<b\le n}\rho_+\left(( a\hspace{0,2cm} b) -\langle a\hspace{0,2cm}b\rangle \right)\\
&=2\rho_-\left((N+1)n+ \sum_{1\le a<b\le n}\left(\la a\, b\ra -(a\, b)\right)\right).\end{align*}
Dividing by $2N$ gives the announced formula.
\end{proof}

\section*{Acknowledgements}
The author is grateful to the anonymous referee for helpful comments that lead to this new organisation of the manuscript.   Many thanks are  due to his former PhD advisor, Thierry Lévy,  for  his insights and for his  work on the preliminary versions of the paper and to Franck Gabriel and Guillaume Cébron for stimulating discussions. We wish also  to thank   Benoît Collins for interesting discussions about this work. 

\def\cprime{$'$} \def\cprime{$'$}

\end{document}